\documentclass{article}
\usepackage{makeidx}
\usepackage{amssymb}
\usepackage{amsmath}

\newtheorem{theorem}{Theorem}
\newtheorem{acknowledgement}[theorem]{Acknowledgement}

\newtheorem{lemma}[theorem]{Lemma}

\newtheorem{proposition}[theorem]{Proposition}

\input{tcilatex}

\begin{document}

\title{Hitting time and dimension in Axiom A systems and generic interval
excanges}
\author{Stefano Galatolo \footnote{Dipartimento di matematica applicata, Universita di Pisa, Italia. email:
galatolo@dm.unipi.it} }
\maketitle

\begin{abstract}
In this note we prove that for equilibrium states of axiom A systems the
time $\tau _{B}(x)$ needed for a typical point $x$ to enter for the first
time in a typical ball $B$ with radius $r$ scales as $\tau _{B}(x)\sim r^{d}$
where $d$ is the local dimension of the invariant measure at the center of
the ball. A similar relation is proved for a full measure set of interval
excanges. Some applications to Birkoff averages of unbounded (and not $L^{1}$
) functions are shown.
\end{abstract}

\section{Introduction}

It is well known by classical recurrence results that a typical orbit of a
dynamical systems comes back (in any reasonable neighborhood) near to its
starting point. The quantitative study of recurrence quantifies the speed of
this coming back, estimating, for example, how much time is needed to come
back in any ball centered in the starting point (the reader can find and
exposition of more and less recent developments about this kind of questions
in the survey \cite{CG}). It turns out that in many cases the scaling law of
return times are related to the dimension of the invariant measure of the
system. More precisely, let us consider a starting point $x$, a ball $B(x,r)$
and the time $\tau _{B(x,r)}(x)$ needed for the starting point $x$ to come
back to $B(x,r).$ With these notations we have, for example, (\cite{S}, \cite%
{STV}) that in exponentially mixing systems or positive entropy (with some
small tecnical assumptions) systems over the interval $\tau _{B(x,r)}(x)\sim
r^{d(x)}$, where $d(x)$ is the local dimension at $x.$ Moreover (\cite{Bo}, 
\cite{BS}) in general finite measure preserving systems the recurrence gives
a lower bound to the dimension ($\underset{r\rightarrow 0}{\lim }\frac{\log
\tau _{B(x,r)}(x)}{-\log r}\leq d(x)$).

A similar and strictly related (see e.g. \cite{LHV}) problem is about the
time needed for a typical point $x$ of an ergodic system to enter in some
neighborhood of another point $y.$This leads to the \emph{hitting time }%
(also called \emph{waiting time})\emph{\ }indicators. For example let us
denote by $\tau _{B(y,r)}(x)$ the time needed for the point $x$ to enter in
the ball $B(y,r)$ with center $y$ and radius $r$ (this in some sense \emph{%
generalizes} the recurrence because we are allowed to consider an arriving
point $y$ different from the starting point $x$). Again we consider the
scaling behavior of this waiting time for small $r$. The waiting time
indicator will have value $R$ if $\tau _{B(y,r)}(x)\sim r^{R}$ (precise
definitions in section 2)$.$ Again, there are relations with the local
dimension. Some general relations are proved in \cite{G} (see thm. \ref{GAN}%
) and show that the waiting time indicator gives an upper bound to the
dimension. Moreover there is a class of systems where the waiting time is
equal to dimension. This class of systems includes for example systems
having exponential distribution of return times in small balls (this
includes many \emph{more or less} hyperbolic systems over the interval \cite%
{BSTV}). We remark that exponential return times in balls is conjecured but
yet not proved in general Axiom A systems, thus equality between hitting
time and dimension for axiom A does not follow from such result.

We want to remark that there are also some relevant cases where the equality
between recurrence or hitting time with dimension does not hold, hence this
kind of questions are not trivial. Such cases includes rotations by
Liuouville numbers (see \cite{BS},\cite{KS}), and Maps having an indifferent
fixed point and infinite invariant measure (\cite{GKP}).

A further motivation for this kind of studies is that the relations between
recurrence (and similar) with dimension are used in the physical literature 
\cite{HJ},\cite{GE},\cite{JKLPS} to provide numerical methods for the study
of the Hausdorff dimension of attractors. Since recurrence gives a lower
bound to dimension and hitting time gives an upper bound, the combined use
of these may produce efficient numerical estimators for the dimension of
attractors.

The main result of this note is to show that in nontrivial nice examples
such as Axiom A systems and typical interval exchanges, the hitting time
indicator equals the local dimension $d(y)$ of the considered measure.
Hence, in such systems, for typical $x$ and $y$ we will have $\tau
_{B(y,r)}(x)\sim r^{d(y)}$. As an application of these results we give an
estimation for the Birkoff averages of functions (in ergodic systems) having
some asymptote and\ no finite $L^{1}$norm. Here the Birkoff average will
increase to infinity as the number of iterations increases (this is
trivially by ergodic theorem). The hitting time indicator will give an
estimation about the speed of going to infinity of such average.

\begin{acknowledgement}
I wish to thank Corinna Ulcigrai , Jean Rene Chazottes and Stefano Marmi for
fruitful discussions, which allowed me to discover some relevant literature
and to simplifiyng much the proof of the main result.
\end{acknowledgement}

\section{Generalities and a criteria for hitting time and dimension}

In the following we will consider a discrete time dynamical system $(X,T)$
were $X$ is a separable metric space equipped with a Borel finite measure $%
\mu $ and $T:X\rightarrow X$ is a measurable map.

Let us consider the first entrance time of the orbit of $x$ in the ball $%
B(y,r)$ with center $y$ and radius $r$ 
\begin{equation*}
\tau _{r}(x,y)=\min (\{n\in \mathbf{N},n>0,T^{n}(x)\in B(y,r)\})\,.
\end{equation*}%
By considering the power law behavior of $\tau _{r}(x,y)$ as $r\rightarrow 0$
let us define the hitting time indicators as%
\begin{equation*}
\overline{R}(x,y)=\mathrel{\mathop{limsup}\limits_{r\rightarrow 0}}\frac{%
\log (\tau _{r}(x,y))}{-\log (r)},\underline{R}(x,y)=\mathrel{%
\mathop{liminf} \limits_{r\rightarrow 0}}\frac{log(\tau _{r}(x,y))}{-\log (r)%
}\,.
\end{equation*}

If for some $r,$ $\tau _{r}(x,y)$ is infinite then $\overline{R}(x,y)$ and $%
\underline{R}(x,y)$ are set to be equal to infinity. The indicators $%
\overline{R}(x)$ and $\underline{R}(x)$ of quantitative recurrence defined
in \cite{BS} are obtained as a special case, $\overline{R}(x)=\overline{R}%
(x,x)$, $\underline{R}(x)=\underline{R}(x,x)$.

We recall some basic properties of $\overline{R}(x,y):$

\begin{proposition}
\label{inizz}$R(x,y)$ satisfies the following properties

\begin{itemize}
\item $\overline{R}(x,y)=\overline{R}(T(x),y)$, $\underline{R}(x,y)=%
\underline{R}(T(x),y)$.

\item If $T$ is $\alpha -Hoelder$, then $\overline{R}(x,y)\geq \alpha 
\overline{R}(x,T(y))$, $\underline{R}(x,y)\geq \alpha \underline{R}(x,T(y))$.

\item If we consider $T^{n}$ instead of $T$ $\ \overline{R}_{T}(x,y)\leq 
\overline{R}_{T^{n}}(x,y)$, $\underline{R}_{T}(x,y)\leq \underline{R}%
_{T^{n}}(x,y)$.
\end{itemize}
\end{proposition}

\emph{Proof. }The first two points comes from \cite{G} (and they comes
directly from definitions). For the third one let us denote with $\tau $ and 
$\tau ^{\prime }$the hitting time with resp. to $T$ and $T^{n}$. By
definition $\overline{R}_{T}(x,y)=\mathrel{\mathop{limsup}\limits_{r%
\rightarrow 0}}\frac{\log (\tau _{r}(x,y))}{-\log (r)}$ but $\tau
_{r}(x,y)\leq n\tau _{r}^{\prime }(x,y)$ and $\frac{\log (\tau _{r}(x,y))}{%
-\log (r)}\leq \frac{\log (\tau _{r}^{\prime }(x,y))+\log n}{-\log (r)},$
and taking the $\lim \sup $ we are done. The same can be done for the $liminf
$.$\square $

In general systems the quantitative recurrence indicator gives only a \emph{%
lower} bound on the dimension (\cite{BS}, \cite{Bo}). The waiting time
indicator instead give an \emph{upper} bound (\cite{G}) to the local
dimension of the measure at the point $y$. This is summarized in the
following

\begin{theorem}
\label{GAN} If $(X,T,\mu )$ is a dynamical system over a separable metric
space, with an invariant measure $\mu .$ For each $y$ 
\begin{equation}
\underline{R}(x,y)\geq \underline{d}_{\mu }(y)\ ,\ \overline{R}(x,y)\geq 
\overline{d}_{\mu }(y)  \label{thm4}
\end{equation}%
holds for $\mu $ almost each $x$. Where $\underline{d}_{\mu }(y)$ and $%
\overline{d}_{\mu }(y)$ are the lower and upper local dimension at $y%
\footnote{%
If $X$ is a metric space and $\mu $ is a measure on $X$ the upper local
dimension at $x\in X$ is defined as $\overline{d}_{\mu }(x)=%
\mathrel{\mathop{limsup}\limits_{r\rightarrow 0}}\frac{log(\mu (B(x,r)))}{%
log(r)}=\mathrel{\mathop{limsup}\limits_{k\in {\bf N}, k\rightarrow \infty}}%
\frac{-log(\mu (B(x,2^{-k})))}{k}$ The lower local dimension $\underline{d}%
_{\mu }(x)$ is defined in an analogous way by replacing $limsup$ with $liminf
$. If $\overline{d}_{\mu }(x)=\underline{d}_{\mu }(x)=d$ almost everywhere
the system is called exact dimensional. In this case many notions of
dimension of a measure will coincide (see for example the book \cite{P}).}.$%
Moreover, if $X$ is a closed subset of ${%
\mathbb{R}
}\mathbf{^{n}}$, then for almost each $x\in X$ 
\begin{equation*}
\overline{R}(x,x)\leq \overline{d}_{\mu }(x)\ ,\ \underline{R}(x,x)\leq 
\underline{d}_{\mu }(x)\,.
\end{equation*}
\end{theorem}

A natural question which is important also from the numerical applications
is whether equality can replace the above inequalities (see the results from 
\cite{S}, \cite{G}, \cite{BS}, \cite{STV} already outlined in the
introduction). The following is a general criteria that assures (together
with theorem \ref{GAN}) for typical points, equality between waiting time
and local dimension.

\begin{lemma}
\label{lemma1} Let $x\in X$ and 
\begin{equation*}
SF_{r}^{n}(x)=X-B(x,r)\cap T^{-1}(X-B(x,r))\cap T^{-2}(X-B(x,r))\cap ...\cap
T^{-n}(X-B(x,r))
\end{equation*}%
be the set of points that after $n$ steps never enters into $B(x,r).$ If for
each $\epsilon >0$ we have $\sum \mu (SF_{2^{-n}}^{\mu
(B(x,2^{-n}))^{-1-\epsilon }})<\infty $ then for almost each $y$ it holds $%
\overline{R}(y,x)\leq \overline{d}_{\mu }(x)$ and \underline{$R$}$(y,x)\leq 
\underline{d}_{\mu }(x).$
\end{lemma}

\emph{Proof. }The proof follows by a Borel Cantelli argument. Let $%
R_{\epsilon }=\{y\in X,\overline{R}(y,x)\geq (1+\epsilon )\overline{d}_{\mu
}(x)\}$. If we prove that this set has measure zero for each $\epsilon $ we
are done. If we know that for some $\epsilon $ $\sum \mu (SF_{2^{-n}}^{\mu
(B(x,2^{-n}))^{-1-\epsilon }})<\infty $ this means that the set of points
such that $\tau _{2^{-n}}(y,x)>\mu (B(x,2^{-n}))^{-1-\epsilon }$ for
infinitely many $n$ has zero measure. Taking logarithms and dividing by $n$
we have $\frac{\log (\tau _{2^{-n}}(x,y))}{n}\leq (1+\epsilon )\frac{\log
(\mu (B(x,2^{-n})))}{-n}$ eventually (as $n$ increases) for a full measure
set and then $\overline{R}(y,x)=\lim \sup \frac{\log (\tau _{2^{-n}}(x,y))}{n%
}\leq (1+\epsilon )\lim \sup \frac{\log (\mu (B(x,2^{-n})))}{-n}=(1+\epsilon
)\overline{d}_{\mu }(x)$ on a full measure set. This is true for each $%
\epsilon $ and we have the statement. The same can be done for the proof of 
\underline{$R$}$(y,x)\leq \underline{d}_{\mu }(x)\square $

\section{Axiom A systems}

In this section we will consider Axiom A systems, we will apply the
properties of Gibbs measures to prove that they satisfy Lemma \ref{lemma1}
at almost all points. We will prove the following

\begin{theorem}
\label{axa} If $X$ is a basic set of an axiom A diffeomorphism$,$ $\mu $ is
an equilibrium measure for an Hoelder potential defined on $X$. Then $%
(X,T,\mu )$ satisfies Lemma \ref{lemma1} at almost each $x\in X$ and hence
for almost each $x$ it holds \underline{$R$}$(y,x)=\underline{d}_{\mu }(x),%
\overline{R}(y,x)=\overline{d}_{\mu }(x)$ for almost each $y.$
\end{theorem}

First we need a general estimation on the behavior of a certain kind of
sequences.

\begin{lemma}
\label{lemma2}Let $0<m<1,$and $a_{n}$ be defined by $\left\{ 
\begin{array}{c}
a_{n}=a_{n-1}m+s_{n} \\ 
a_{0}=m^{2}%
\end{array}%
\QQfntext{-1}{
Keywords: dimension, quantitative recurrence, Axiom A, interval exchanges,
Birkoff sums}\right. $ where $s_{n}=\frac{2n+1}{n^{2}(n+1)^{2}}=\frac{1}{%
i^{2}}-\frac{1}{(i+1)^{2}}$ then for $n\geq 2$ it holds $a_{n}\leq \frac{m^{[%
\frac{n}{2}]}}{1-m}+\frac{4}{n^{2}}.$
\end{lemma}

\emph{Proof. }We have

\begin{equation*}
a_{n}=m^{n+1}+m^{n-1}s_{1}+m^{n-2}s_{2}+m^{n-3}s_{3}+...+ms_{n-1}+s_{n}.
\end{equation*}

Since $s_{i}<1$ and $m<1$ then $a_{n}\leq
\sum_{n/2}^{n}m^{i}+\sum_{n/2}^{n}s_{i}=\frac{m^{[\frac{n}{2}]}-m^{n}}{1-m}+%
\frac{1}{([n/2])^{2}}-\frac{1}{([n/2]+1)^{2}}\leq \frac{m^{[\frac{n}{2}]}}{%
1-m}+\frac{4}{n^{2}}.\square $

\emph{Proof of Theorem \ref{axa}. }We already know that (thm. \ref{GAN}) $%
\underline{R}(x,y)\geq \underline{d}_{\mu }(y)\ ,\ \overline{R}(x,y)\geq 
\overline{d}_{\mu }(y)$. For the opposite inequalities, first we remark that
(see \cite{Bow} pag. 72) $X=X_{1}\cup ...\cup X_{l}$ where $%
T(X_{i})=X_{i+1},\ (T(X_{l})=X_{1})$ and $T^{l}|_{X_{i}}$ is topologically
mixing.

By Lemma \ref{inizz} we can suppose that $x,y$ belongs to the same $X_{i}.$
and that $T$ is topologically mixing (replacing $T$ with $T^{l}$ we have a
mixing transformation on $X_{i}$, moreover by the fourth point of Lemma \ref%
{inizz} we see that if we have an upper bound for $\overline{R}_{T^{l}}$ and 
\underline{$R$}$_{T^{l}}$then this is also an upper bound for $\overline{R}%
_{T}(x,y)$ and \underline{$R$}$_{T}(x,y).$ Since in this proof we are
looking for an upper bound, by replacing $T$ with $T^{l}$ we can suppose
that the map is topologically mixing).

To estimate the measure of the set $SF_{r}^{n}(x)$ let us consider a Markov
partition $Z=\{Z_{i}\}$ of $X.$ Let $Z_{m}^{n}=T^{-m}(Z)\vee ...\vee
T^{-n}(Z).$ By uniform hyperbolicity there are constants $C,\lambda >0$ such
that $\ diam(Z_{-n}^{n})\leq Ce^{-\lambda n}.$ By this \ we know that there
is some%
\begin{equation*}
n(r)\leq -\lambda ^{-1}(\log r-\log C)
\end{equation*}%
such that the partition is of size so small that there is one element $Z_{0}$
of the partition $\overline{Z}=Z_{-n(r)}^{n(r)}$ which is included in $%
B(x,r).$

Now $SF_{r}^{n}(x)\subset B_{0}^{n}=X-Z_{0}\cap T^{-1}(X-Z_{0})\cap ...\cap
T^{-n}(X-Z_{0}).$ We remark that $B_{0}^{n}$ is the union of many cylinders.
The measure of $B_{0}^{n}$ decreases very fast by the weak Bernoulli
property of the equilibrium measure $\mu .$ Indeed by \cite{Bow} pag. 90, we
know that since the map $T$ can be supposed to be topologically mixing and
then $\mu $ has the weak Bernoulli property: i.e. let us consider $t,s\geq 0$
and the partitions $P_{s}=Z\vee T^{-1}(Z)\vee ...\vee T^{-s}(Z)$ and $%
Q_{t}=T^{-t}(Z)\vee ...\vee T^{-t-k}(Z).$ For each $\epsilon ,$ if $%
t-s=N_{Z}(\epsilon )$ is big enough, then $\underset{P\in P_{s},Q\in Q_{t}}{%
\sum }\mu (Q\cap P)-\mu (P)\mu (Q)<\epsilon .$ Moreover by \cite{Bow},
theorem 1.25 we can find an estimation for $N_{Z}(\epsilon )$ as a function
of $\epsilon $ (see \cite{Bow} pag. 38): $N_{Z}(\epsilon )=-c\log (\epsilon
)+c^{\prime },$ where $c,c^{\prime }$ are constants depending on $\mu ,T$
and $Z$.

The estimation for $\mu (B_{0}^{l})$ follows from the fact that a non empty
cylinder for the partition $\overline{Z}=Z_{-n(r)}^{n(r)}$ is also a
cylinder for the partition $Z.$

Indeed the cylinder $\overline{\mathbf{z}}_{m}=\overline{Z}_{i_{1}}\cap
T^{-1}(\overline{Z}_{i_{2}})\cap ...\cap T^{-m-1}(\overline{Z}_{i_{m}})$ , $%
\overline{Z_{i}}\in \overline{Z}$ satisfies $\overline{\mathbf{z}}_{m}=%
\mathbf{z}_{m+2n}$ where $\mathbf{z}_{m+2n}=T^{n}(Z_{j_{1}})\cap
T^{n-1}(Z_{j_{2}})\cap ...\cap T^{-m-1-n}(Z_{j_{(m+2n)}})$ is a cylinder of $%
Z$ and $\overline{Z}_{i_{k}}=T^{-k+n}Z_{j_{k}}\cap ...\cap
T^{-k-n}Z_{j_{k+2n}}.$ Since $\mu $ is preserved then $\mu
(T^{-k+n}Z_{j_{k}}\cap ...\cap T^{-k-n}Z_{j_{k+2n}})=\mu
(T^{-k}Z_{j_{k}}\cap ...\cap T^{-k-2n}Z_{j_{k+2n}})$ hence we can apply the
weak Bernoulli property to such cylinders obtaining that also $\overline{Z}$
satisfies such a property and 
\begin{equation}
N_{\overline{Z}}(\epsilon )\leq c\log (\epsilon )+c^{\prime }+2n
\label{eqnr}
\end{equation}%
(where $c,c^{\prime }$ are the \ constants of $N_{Z}(\epsilon )$ as above).
\ We recall that $n$ depends on $r$ and we can choose $n(r)\leq -\lambda
^{-1}(\log r-\log C).$

Finally let us apply the weak Bernoulli property of $\overline{Z}$ to get an
estimation for $\mu (B_{0}^{n}).$ Let us set $m=\mu (X-Z_{0})$ and $\epsilon
(i)=\frac{2i+1}{i^{2}(i+1)^{2}}=\frac{1}{i^{2}}-\frac{1}{(i+1)^{2}}$. We
have (eq. \ref{eqnr}) that setting $C^{\prime }(r)=c^{\prime }-2\lambda
^{-1}(\log r-\log C)$ then $N_{\overline{Z}}(\epsilon )\leq c\log (\epsilon
)+c^{\prime }+2n=c\log (\frac{2i+1}{i^{2}(i+1)^{2}})+C^{\prime }(r)$ and
there is a $C$ s.t. $N_{\overline{Z}}(\epsilon )\leq C\log (i)+C^{\prime
}(r).$

Let us set $n_{i}=\sum_{j\leq i}N_{\overline{Z}}(\epsilon (j)).$ Thus for
each $\delta $ there is a $K$ such that$,$ if $i$ is big enough 
\begin{equation}
n_{i}\leq -Ki^{1+\delta }\log r.  \label{eq2}
\end{equation}%
The measure of $B_{0}^{n_{i}}$ can then be estimated applying $i$ times the
Bernoulli property, with $\epsilon (i)=$ $\frac{2i+1}{i^{2}(i+1)^{2}}$ as
above, to \ subcylinders of increasing length $N_{\overline{Z}}(\epsilon
(1)),N_{\overline{Z}}(\epsilon (1))+N_{\overline{Z}}(\epsilon (2)),N_{%
\overline{Z}}(\epsilon (1))+N_{\overline{Z}}(\epsilon (2))+N_{\overline{Z}%
}(\epsilon (3))...$ obtaining by the Bernoulli property of $\mu $

\begin{eqnarray*}
\mu (B_{0}^{N(\epsilon (1))}) &\leq &m^{2}+\epsilon (1), \\
\mu (B^{N(\epsilon (1))+N(\epsilon (2))}) &\leq &(m^{2}+\epsilon
(1))m+\epsilon (2), \\
\mu (B^{N(\epsilon (1))+N(\epsilon (2))+N(\epsilon (3))}) &\leq
&((m^{2}+\epsilon (1))m+\epsilon (2))m+\epsilon (3),...\ 
\end{eqnarray*}

$\ $Hence \ by Lemma \ref{lemma2} above

\begin{equation*}
\mu (B_{0}^{n_{i}})\leq \frac{m^{[\frac{i}{2}]}}{1-m}+\frac{4}{i^{2}}.
\end{equation*}

We remarked that $SF_{r}^{n}\subset B_{0}^{n}.$ If we consider another
element $Z_{1}$ of $\overline{Z}$ with $Z_{1}\subset B(x,r)$ and $B_{1}^{n}=$
$X\cap T^{-1}(X-(Z_{0}\cup Z_{1}))\cap ...\cap T^{-n}(X-Z_{0}\cup Z_{1})),$
we have also $SF_{r}^{n}\subset B_{1}^{i}\subset B_{0}^{l}.$ Now considering
a sequence $Z_{0},...,Z_{w}$ of elements of $\overline{Z}$ with $%
Z_{0},...,Z_{w}\subset B(x,r)$ and $B_{w}^{n}=X\cap T^{-1}(X-(Z_{0}\cup
...\cup Z_{w}))\cap ...\cap T^{-n}(X-(Z_{0}\cup ...\cup Z_{w})),$ we have
also $SF_{r}^{n}\subset B_{w}^{n}.$ The measure of $B_{w}^{n}$ can be
estimated as above, obtaining $\mu (B_{w}^{n_{i}})\leq \frac{m_{w}^{i/2}}{%
1-m_{w}}+\frac{4}{i^{2}},$ where $m_{w}=\mu (X-(Z_{0}\cup ...\cup Z_{w})).$

Now, refining again the partition $\overline{Z}$ if necessary (this is true
because the diameter of each $Z_{i}$ is less or equal than $Ce^{-\lambda n},$
and this will only change the constants in $N(\epsilon )$) we can suppose
that the diameter of each piece of the partition has diameter less than $%
r/4. $ We then have that we can choose $Z_{0},...,Z_{w}$ such that $B(x,%
\frac{r}{2})\subset Z_{0}\cup ...\cup Z_{w}\subset B(x,r).$ Then $\mu
(X-(Z_{0}\cup ...\cup Z_{w}))\leq \mu (X-B(x,\frac{r}{2})).$This gives,%
\begin{equation*}
\mu (B_{w}^{n_{i}})\leq \frac{(1-\mu (B(x,\frac{r}{2})))^{\frac{i}{2}}}{\mu
(B(x,\frac{r}{2}))}+\frac{4}{i^{2}}.
\end{equation*}

By \ref{eq2} we then have $i\geq \frac{n_{i}^{\frac{1}{1+\delta }}}{(K\log
r/4)^{^{\frac{1}{1+\delta }}}}$, by this

\begin{eqnarray*}
\mu (SF_{2^{-n}}^{\mu (B(x,2^{-n}))^{-1-\epsilon }}) &\leq &\mu (B_{w}^{\mu
(B(x,2^{-n}))^{-1-\epsilon }})\leq  \\
&\leq &\frac{(1-\mu (B(x,2^{-n-1})))^{(Kn+\log 4)^{\frac{-1}{1+\delta }}\mu
(B(x,2^{-n}))^{\frac{-1-\epsilon }{1+\delta }}}}{\mu (B(x,2^{-n-1}))}+ \\
&&+\frac{4}{(Kn+\log 4)^{\frac{-2}{1+\delta }}\mu (B(x,2^{-n-1}))^{^{\frac{%
-2-2\epsilon }{1+\delta }}}}.
\end{eqnarray*}%
Since in our case $d_{\mu }(x)=d<\infty \ a.e.,$ there is a constant $Q,$
(which depends on the local dimension) such that $0<Q<\frac{\mu
(B(x,2^{-n-1}))}{\mu (B(x,2^{-n}))}<1$ when $n$ is big, recalling that $%
\delta $ can be chosen as small as we want and hence smaller than $\epsilon $%
, then $\mu (SF_{2^{-n}}^{\mu (B(x,2^{-n}))^{-1-\epsilon }})$ is less than
about $\frac{(e^{-\frac{Q}{2}})^{Kn^{\frac{-1}{1+\delta }}\mu (B(x,2^{-n}))^{%
\frac{\epsilon -\delta }{1+\delta }}}}{\mu (B(x,2^{-n-1}))}+4(Kn+\log 4)^{%
\frac{2}{1+\delta }}\mu (B(x,2^{-n-1}))^{^{\frac{2+2\epsilon }{1+\delta }}}$
and we have $\underset{n}{\sum }\mu (SF_{2^{-n}}^{\mu
(B(x,2^{-n}))^{-1-\epsilon }})<\infty .$This is enough to apply the Lemma %
\ref{lemma1} and have the required statement.$\square $

\section{Interval exchanges}

An interval excange is a piecewise isometry which preserves the Lesbegue
measure. In this section we apply a \ result of Boshernitzan about a full
measure class of uniquely ergodic interval exchanges maps to prove equality
between hitting time and dimension at discontinuity points. We refer to \cite%
{Bo2} for generalities on this important class of maps. 

\begin{theorem}
\label{IET}For a typical interval exchange transformation $\ T$ (for a full
measure set, in the space of interval exchanges) for each discontinuity
point $x_{0}$ it holds \underline{$R$}$(y,x_{0})=1$ for almost each $y\in
\lbrack 0,1].$
\end{theorem}

\emph{Proof. }By a result of (\cite{Bo2}) we have that if $T$ is a typical
i.e.t. and $\delta (n)$ is the minimum distance between the discontinuity
points of $T^{-n},$ then there is a\ constant $C$ and a sequence $n_{k}$
such that $\delta (n_{k})\geq \frac{C}{n_{k}}.$ If $x_{0}$ is a
discontinuity point then for each $n_{k}$ it also hold that $\min_{i,j\leq
n_{k}}d(T^{-i}(x_{0}),T^{-j}(x_{0}))\geq \frac{C}{n_{k}}.$ Let us consider
the set $J_{k}=\cup _{i\leq n_{k}}B(T^{-i}(x_{0}),\frac{C}{3n_{k}})$. Since
it is a union of disjoint balls the measure of $J_{k}$ is $\frac{2C}{3}.$
This implies that in the interval $[0,1]$ there is a positive measure set $J$
of points belonging to infinitely many $J_{k}.$

If a point $y$ is in $J$ then for a subsequence $n_{k_{i}}$ it holds $%
d(T^{n_{k_{i}}}(y),x_{0})\leq \frac{C}{3n_{k_{i}}}$. This is true because by
the Boshernitzan result there are no counter images of other discontinuity
points in the interval $[y,T^{-n}(x_{0})]$ and then these two points cannot
be separated during $n$ iterations of the map.

Since $\lim \inf \ n\ d(T^{n}(y),x_{0})\leq \frac{C}{3}$ then \underline{$R$}%
$(y,x_{0})\leq 1$ and we have the statement for $y$ varying in a positive
measure set $J$. Since the system is ergodic, by Proposition \ref{inizz} we
have the statement.$\square $

Since in interval exchanges th only source of initial condition sensitivity
is the discontinuity (the orbits of two points can be only separated by a
discontinuity) we remark that an estimation of the approaching speed of
typical orbits to the discontinuity is useful to estimate the kind of "weak"
chaos that is present is such maps. The theorem above in some sense can give
(using the construcion done in \cite{BGI} ) an upper bound on the initial
condition sensitivity of such maps. We will not go into details about this
in this work however.

\section{Hitting time and Birkoff sums}

Let us consider a function $f:X-\{x_{0}\}\rightarrow 
\mathbb{R}
$ , $f\geq 0$ which is continuos and which satisfies $\int_{X}fd\mu =\infty $
because it has an asymptote in $x_{0}$ where $f(x)\sim d(x,x_{0})^{-\alpha
}. $

By the ergodic theorem we know that for almost each $x$ the Birkoff average $%
\frac{S_{n}(x)}{n}=\frac{1}{n}\sum_{i=0}^{n}f(T^{i}(x))$ is such that $\frac{%
S_{n}(x)}{n}\rightarrow \infty .$ By the results of the previous sections,
if we know the hitting time indicator at  $x_{0}$ we can have an estimation
for the speed of increasing of $\frac{S_{n}(x)}{n}.$

\begin{theorem}
Let us suppose that near $x_{0}$ we have $\ 0<\lim \frac{f(x)}{%
d(x,x_{0})^{-\alpha }}<\infty $, for $\alpha >1$ then for almost each $x$

\begin{equation*}
\frac{\alpha }{\overline{R}(x,x_{0})}\leq \underset{n\rightarrow \infty }{%
\lim \sup }\frac{\log (S_{n}(x))}{\log (n)}\leq \frac{\alpha }{\underline{R}%
(x,x_{0})}+1
\end{equation*}
\end{theorem}

\emph{Proof. }By the definition of \underline{$R$}$(y,x_{0})$ we obtain (see 
\cite{G}) that for each $\epsilon >0$ 
\begin{equation*}
\underset{n\rightarrow \infty }{\lim \inf }n^{\frac{1}{\underline{R}(y,x_{0})%
}+\epsilon }d(T^{n}(y),x_{0})=\infty .
\end{equation*}%
Then we have that if $n$ is big enough $d(T^{n}(y),x_{0})\geq n^{-\frac{1}{%
\underline{R}(y,x_{0})}-\epsilon }.$ Now we remark that since $\ X$ is
compact there are $c_{1},c_{2}$ such that $f(x)\leq \max
(c_{1},c_{2}d(x,x_{0})^{-\alpha }).$ Then if $n$ is big enough

\begin{eqnarray*}
\sum_{i=0}^{n}f(T^{i}(y)) &\leq &\sum_{i=0}^{n}\max
(c_{1},c_{2}d(T^{i}(y)-x_{0})^{-\alpha })\leq  \\
&\leq &\sum_{i=0}^{n}\max (c_{1},c_{2}n^{\frac{\alpha }{\underline{R}%
(y,x_{0})}+\alpha \epsilon })\leq nc_{1}+c_{2}n^{\frac{\alpha }{\underline{R}%
(y,x_{0})}+\alpha \epsilon +1}
\end{eqnarray*}%
and we have $\underset{n\rightarrow \infty }{\lim \sup }\frac{\log (S_{n}(x))%
}{\log (n)}\leq \frac{\alpha }{\underline{R}(y,x_{0})}+1.$ On the other
hand, by the definition of $\overline{R}(y,x_{0})$ we have that frequently $%
d(T^{n}(y),x_{0})\leq n^{\frac{-1}{\overline{R}(y,x_{0})}+\epsilon }$, then
frequently $\sum_{i=0}^{n}f(T^{i}(y))\geq cn^{\frac{\alpha }{\overline{R}%
(y,x_{0})}-\alpha \epsilon }$.$\square $

By the above result and the previous ones it easily follows :

\begin{enumerate}
\item (by thm. \ref{GAN}) In a general system, if the local dimension at $%
x_{0}$ is $d_{\mu }(x_{0}).$ Then for almost each $x$ 
\begin{equation*}
\underset{n\rightarrow \infty }{\lim \sup }\frac{\log (S_{n}(x))}{\log (n)}%
\leq \frac{\alpha }{d_{\mu }(x_{0})}+1
\end{equation*}

\item (by thm. \ref{IET}) If $T$ is an IET \ and $x_{0}$ is a discontinuity
point then 
\begin{equation*}
\ \alpha \leq \underset{n\rightarrow \infty }{\lim \sup }\frac{\log
(S_{n}(x))}{\log (n)}\leq \alpha +1
\end{equation*}

\item (by thm. \ref{axa}) If $(X,T)$ is axiom A (with an equilibrium
measure, as in thm. \ref{axa}), $x_{0}$, $x$ are typical and $d$ is the
dimension of the measure then 
\begin{equation*}
\frac{\alpha }{d}\leq \underset{n\rightarrow \infty }{\lim \sup }\frac{\log
(S_{n}(x))}{\log (n)}\leq \frac{\alpha }{d}+1.
\end{equation*}
\end{enumerate}

\end{document}